\documentclass[AMA,STIX1COL]{WileyNJD-v2}

\usepackage{hyperref}

\articletype{RESEARCH ARTICLE}%

\received{26 April 2016}
\revised{6 June 2016}
\accepted{6 June 2016}

\newtheorem{defi}{Definition}
\newtheorem{teo}{Theorem}

\usepackage{epsfig}

\begin{document}

\title{New quality measures for quadrilaterals and new discrete functionals for grid generation}

\author[1]{Gonz\'alez Flores, Guilmer*}

\author[2]{Barrera S\'anchez, Pablo}

\authormark{Gonz\'alez Flores, Guilmer \textsc{et al}}

\address[1]{\orgdiv{Facultad de Ciencias}, \orgname{Universidad Nacional Aut\'onoma de M\'exico}, \orgaddress{\state{CDMX}, \country{M\'exico}}}

\corres{Guilmer Gonz\'alez Flores,     Laboratorio de C\'omputo Cient\'\i fico, 
    Departamento de Matem\'aticas, Facultad de Ciencias, 
    Universidad Nacional Aut\'onoma de M\'exico, 
    Ciudad universitaria, 04510 Ciudad de M\'exico, M\'exico. \\ \email{guilmerg@ciencias.unam.mx} and \\  \email{pablobarrera@ciencias.unam.mx} }

%\presentaddress{This is sample for present address text this is sample for present address text}

\abstract[Summary]{In this paper, we review some grid quality metrics \cite{Robinson1987, Lo1989, Field2000, Knupp2001, Remacle2012} and define some new quality measures for quadrilateral elements. Usually, a maximum value of  a quality measure corresponds to the minimum value of the energy density over the grid \cite{Ivanenko2000}. 

We also define new discrete functionals which are implemented as objective functions in an optimization-based method for quadrilateral grid generation and improvement. These  functionals are linearly combined with a discrete functional whose domain has an infinite barrier at the boundary of the set of unfolded grids like $S_{\omega,\epsilon}(G)$, see \cite{Barrera2010}, in order to preserve convex grid cells in each step of the optimization process. }


\jnlcitation{\cname{%
\author{G. Gonz\'alez Flores}, and
\author{P. Barrera S\'anchez}}  (\cyear{2022}), 
\ctitle{New quality measures for quadrilaterals and new discrete functionals for grid generation}, \cjournal{xxxxx.}, \cvol{xxxx}.}

\maketitle

\footnotetext{\textbf{Abbreviations:} ANA, anti-nuclear antibodies; APC, antigen-presenting cells; IRF, interferon regulatory factor}

\section{Introduction}
There is a relatively large number of papers on numerical grid generation of unstructued grids. There is also a special interest in studying  meshes formed by triangular elements. Our interest here is to generate structured meshes with quadrilateral elements, however, all the disscussion can be applied to unstructured meshes.

The simplest way to generate a structured  mesh is by interpolation of the boundaries, but it is difficult to ensure that the mesh thus obtained is a convex one. In \cite{Barrera2010}, the authors made a review of some functionals and conditons that guarantee the existence of optimal meshes which are $\epsilon$-convex over irregular planar regions. Our interest now is to improve mesh quality via the control the shape of the elements. The improvement of mesh quality can be done in two ways:

\begin{description}
\item [{\em Clean-up}.] Basically it consists in elimination, insertion and  reconnection of nodes in order to eliminate the worst elements. Some authors call that this procedure as topological optimization in the sense that the connectivity of the nodes is removed to obtain an optimal configuration.
\item [{\em Smoothing}.] It consists on node reposition without changing the connectivity of the elements. 
\end{description}

In both cases, the goal is to obtain a quality mesh with a low number of distorted elements. To achieve this goal for a quadrilateral it is neccesary to define an {\em ad hoc} quality measure.

\begin{defi}
We say that a real valued function $\mu(Q)$ over a quadrilateral $Q$ is a quality measure in the sense of Field-Oddy, if it

\begin{enumerate}
\item[1)] has the ability to detect degenerate elements;
\item[2)] is bounded and continuous;
\item[3)] is independent of scale;
\item[4)] is normalized;
\item[5)] and it is invariant under rigid transformations;
\end{enumerate}

\end{defi}
For practical purposes it is convenient to define an acceptability interval $[\mu_0, 1]$ for the quality value when a quadrilateral has a suitable shape. When the quality value lies outside this interval the quadrilateral is distorsioned.

In this work we are interested in identifying the shape of the cells and quantifying the distortion of a quadrilateral when it is not a square or a rectangle. For this purpose, we will review the most used quality measures for rectangles and then we will propose new quality measures.

\section{Background}
Following the ideas behind the quality measures reported for triangles, it is straightforward to define some figures that measure the shape of the quadrilaterals. One of these is the aspect ratio, which is defined  by comparing with the ideal case when the quadrilateral is a rectangle: the ratio of the largest to the smallest sides. An estimator for this ratio was discussed in 1987 by Robinson \cite{Robinson1987}. The idea is to associate a rectangle to the  convex quadrilateral: a rectangle passing through the midpoints of the sides of the quadrilateral, see the figure bellow.

\begin{figure}[hbt]
\centerline{\scalebox{.3}{\includegraphics{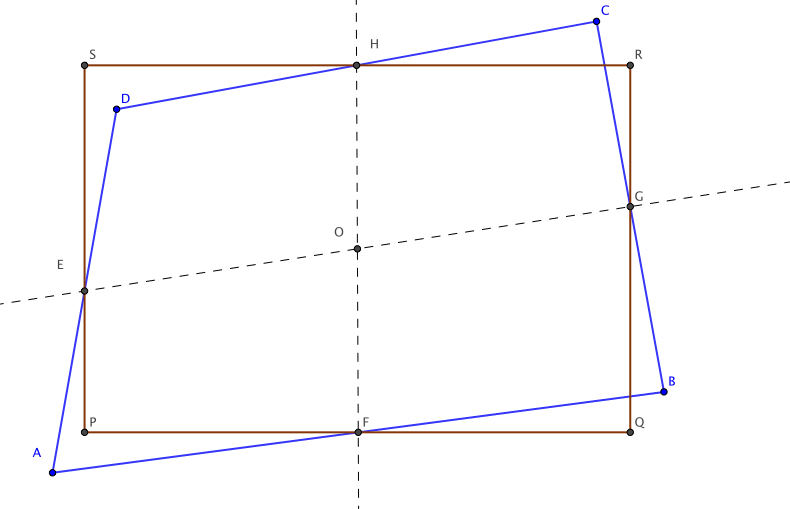}}}
\caption{A rectangle associated with a quadrilateral.}
\end{figure}

This idea is usual in continuum mechanics. Robinson proposes a practical way of calculating it by means of the bilinear mapping between the unit square and the quadrilateral

\begin{eqnarray}
x & = & e_1 + e_2\xi + e_3\eta + e_4 \xi\eta, \\
y & = & f_1 + f_2\xi + f_3\eta + f_4 \xi\eta,
\end{eqnarray}
which yields

$$
\mbox{aspect ratio}  =  \max \left\{ \frac{e_2}{f_3}, \frac{f_3}{e_2} \right\} .
$$

The associated rectangle has sides which are parallel to the coordinate axes and pass through the midpoints of the sides of the quadrilateral. In spite of its simplicity, this analytic representation is not satisfactory since it depends on a orthogonal coordinate system. In 2000, Field \cite{Field2000} reviewed this definition and suggested to calculate the aspect ratio of Robinson orthogonalizing the main axes and proposes a quality measure to detect squares.

In 1989, Lo \cite{Lo1989} reviewed the quality measure for triangules $T(a,b,c)$ with side lengths $l_1$, $l_2$ and $l_3$

$$
g_i = 4\sqrt{3}\frac{\mbox{\'area}(T_i)}{l_1^2+l_2^2+l_3^2}
$$
which attains its optimum value in equilateral triangles and proposed again to calculate each one of those values over the four $T_i$ triangles  which are defined by the sides and diagonals of a quadrilateral, but reordering these quantitites in such a way that 

$$
g_1 \leq g_2 \leq g_3 \leq g_4
$$
and use

$$
\mu(Q)  = \frac{g_1 g_2 }{g_3 g_4}.
$$
as quality measure. Lo points out that the optimal value of 1 is obtained for rectangles. This is a quality measure because it is continuous, bounded and identifies degenerate and even non-convex quadrilaterals. The measure that Lo uses for triangles is the reciprocal of the number of condition of a linear mapping $\mu(Q)=1/\kappa_2(Q)$, that Knupp \cite{Knupp2001} used in 2001 to measure the distortion of the elements. Locally, Lo's measure may have more critical points which can be far from representing a rectangle.

Another measure of quality, described by van Rens {\em et al.} \cite{vanRens1998} for quadrilaterals, is given as follows: compute the inner angles $\theta_k$ and define

$$
\mu(Q) = \prod_{k=1}^4 \left (1 - \left| \frac{\frac{\pi}{2}-\theta_k}{\frac{\pi}{2}} \right| \right ) .
$$
This function is continuous, dimensionless and $0\leq \mu(Q) \leq 1$. One can see that $\mu(Q)=0$ if $Q$ is a triangle and  $\mu(Q)=1$ only if $Q$ is a rectangle.

In 2012, Remacle {\em et al.} \cite{Remacle2012} described the Blossom-Quad algorithm to construct a non-structured  mesh with quadrilaterals elements obtained from a previous triangulation and used a cost-function to produce a quality mesh. They used

$$
\mu(Q) = \max \left\{1-\frac{2}{\pi}\max_k \left\{ |\frac{\pi}{2} - \theta_k|\right\} , 0 \right\}
$$
and observed that the value of this function is 1 the $Q$ is a perfect quadrilateral and 0 if any of the angles is greater than or equal to $\pi$, {\em id est} when the quadrilateral degenerates into a triangle or is nonconvex. This function is also a quality measure.

As noted, unlike other measures for rectangles we have discussed up to this point, the two last ones do not depend neither on the shape of the quadrilaterals, nor the aspect ratio or proportion of their sides; they only measure how near or far away is a quadrilateral from being a rectangle using only the internal angles.

Another function based on inner angles was proposed by Wu \cite{Wu2011}. This author used the same idea of Lo: to order the inner angles $\theta_i$ so that $\theta_1 \leq \theta_2 \leq \theta_3 \leq \theta_4$ and define

$$
\mu(Q)  = \frac{\theta_1 \theta_2 }{\theta_3 \theta_4}.
$$
Wu observed that this function attains its optimal value of 1 on rectangles. However, this is not a good measure in the sense of Field-Oddy, since it is not capable of detecting degenerate quadrilaterals.

\section{New quality measures}
In the previous section, we have reviewed some measures that characterize rectangles and we have also pointed out some intervals of acceptability to decide whether that a quadrilateral is close to having the desired shape is distorted. However, we have not noticed how it is that square or rectangle.

\subsection{Minrect 2015}
A very interesting problem in computational geometry is the following: given a cloud of points, calculate the rectangle of the minimum area that contains them. It is known that this problem can be raised directly on the convex hull of the cloud of points and therefore the problem can be regarded as to calculate a rectangle of minimum area that contains a polygon convex.

\begin{figure}[hbt]
\centerline{\begin{tabular}{cc}
\scalebox{.25}{\includegraphics{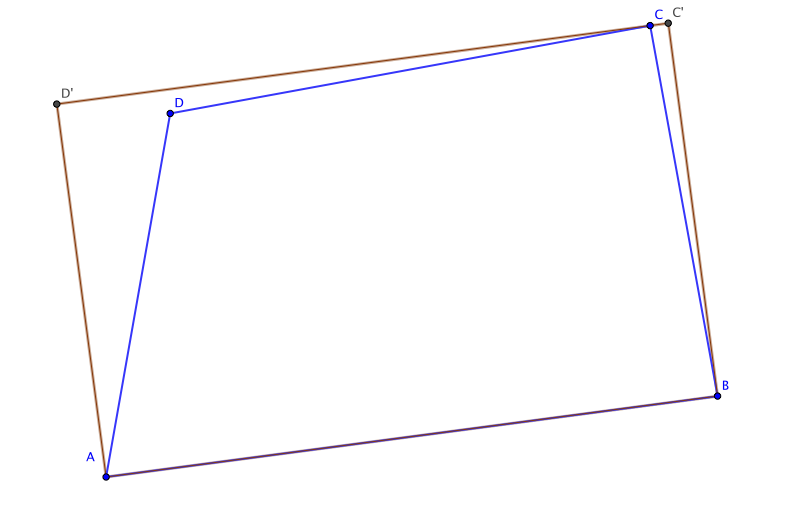}} &
\includegraphics[scale=.25]{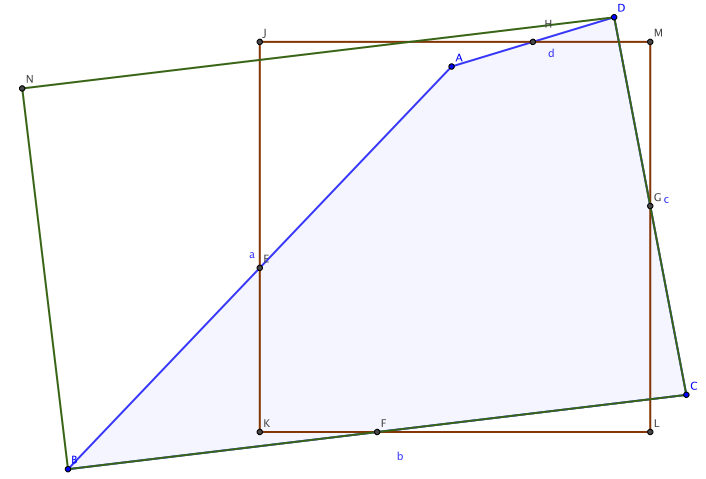} \\
(a) & (b)
\end{tabular}}
\caption{(a) The rectangle of minimum area for one quadrilateral. (b) Comparaci\'on entre {\em aspect} ratio of Robinson and the {\em aspect ratio} of rectangle of minimum area}
\label{minrect}
\end{figure}

We propose to use the rectangle of minimum area to define a distortion measure of the quadrilateral in the sense that measures how close or far is a quadrilateral $Q$ of being a rectangle (see \ref{minrect}(a)).

On one hand,  the cell area $a_c$,  is less than the rectangle area; that is, $a_c\leq a_R$,  and is easy to see that

$$
\frac{2a_c - a_R}{a_R} \leq 1;
$$
see Lassak \cite{Lassak1993} for the proof. The quotient thus defined reaches its maximum value of 1 on rectangles.  Therefore,  we propose the value
 
$$
\mu(Q) =\frac{2a_c - a_R}{a_R}
$$
as a quality measure to characterize rectangles.

One must note that $\mu(Q)$ is a good quality measure according to Field-Oddy,  and $\mu(Q)=0$ when $Q$ is a triangle. 

Another good quality measure in this sense is
$$
\mu(Q) =\frac{2a_{-}}{a_R}
$$
where
\begin{equation}
a_{-} =\min \left\{a_1,\, a_2,\, a_3,\, a_4 \right\}
\label{amin}
\end{equation}
is the minimum area between the four triangles defind by taking the four vertices of a quadrilateral into groups of three. 

\subsection{New aspect ratio}
Using the rectangle of minimum area for $Q$, we propose to use the ratio of the largest to the smallest side as aspect ratio. This measures is invariant under rigid and scaling transformations.

This measure is better than the Robinson aspect ratio. It is easy to construct an example for which the Robinson {\em aspect ratio} 1 but very distorted

{\bf Example :} For the quadrilateral  $A(3.53, 10.21),  B(-10, -4)$, $C(11.81, -1.38)$ y $D(9.27, 11.94)$ the Robinson's {\em aspect ratio} is $1.00$ and using the rectangle of minimum area our {\em aspect ratio} is $1.62$.  See the figure~\ref{minrect}(b).

\subsection{Rectangles 2015}
As we have discussed, some measures to characterize rectangles are based on the inner angles. Another way to do this is to ask $Q$ to be a parallelogram and one of its inner angles to be a right one. This is, we use a measure that imposes a particular condition of a rectangle instead of one on the form of $Q$ .

Our interest is to characterize the rectangles geometrically. A  well  known result in the literature is as follows:

\begin{teo} \label{teoa}
Let $Q$ be a quadrilateral of vertices $A, B, C$ and $D$ whose sides are $a, b, c$ and $d$. The quadrilateral $Q$ is a rectangle if and only if the area of the quadrilateral is written as

$$a_R = \frac{1}{2}\sqrt{(a^2+c^2)(b^2+d^2)}$$
\end{teo}
The proof of this result can be found in Josefsson \cite{Josefsson2013}. The interesting fact about this theorem is that it provides of an analytical expression of the area of a hypothetical rectangle formed by the sum of the square of the opposite sides of $Q$ and compare the square of the area of $Q$ to identify how far if is from being a rectangle.

On the other hand, following the proof of the theorem, it is easy to see that the area $a_c$  of any convex quadrilateral satisfies

$$
2 a_c \leq \sqrt{(a^2+c^2)(b^2+d^2)}.
$$
Using this idea we propose the measure

$$
\mu(Q)  = \frac{2a_{-}}{\sqrt{(a^2+c^2)(b^2+d^2)}}
$$
where $a_{-}$ is defined in equation (\ref{amin}).

This is a good quality measure in the Field-Oddy sense, since it is continuous,  bounded, and capable of indentifying degenerate quadrilaterals (to triangle) as well as to identify if a quadrilateral is non-convex. This measure reaches its optimal value of 1 for rectangles.

An acceptability interval to consider that the quadrilateral is rectangle under this measure is $[.95,\, 1]$.

\subsection{Cuadrados:  Harmonic mean 2017}
An ideal mesh is one in which its cells are close to being squares. If $\mu(T)$ is a good quality measure for triangles, the harmoic mean of the four triangles $T_i$ 

\begin{equation}
\mu(Q) =  \sigma \frac{4}{\sum_{i=1}^4 \frac{1}{\mu(T_i)}}
\label{marmonica}
\end{equation}
is a good quality measure for quadrilaterals, because it is continuous, bounded, invariant under rigid and scaling transformations and identifies degenerate and even non-convex quadrilaterals
since it inherits those properties from $\mu(T)$. Here $\sigma$ is a normalization parameter. 

Rewriting (\ref{marmonica}) we obtain
\begin{equation}
\mu(Q) = \frac{4\sigma \mu(T_1) \mu(T_2) \mu(T_3) \mu(T_4)}{\mu(T_2) \mu(T_3) \mu(T_4)+\mu(T_1) \mu(T_3) \mu(T_4)+\mu(T_1) \mu(T_2)  \mu(T_4)+\mu(T_1) \mu(T_2) \mu(T_3) }
\label{marmonica2v}
\end{equation}
To characterize squares we require a property $\mu(T)$ as seen from

\begin{teo}
If $\mu(T)$ is a good quality measure for triangles according to Field-Oddy in which for isosceles triangle the highest energy among all right triangles is achieved, sothe harmoic mean (\ref{marmonica2v}), characterize squares at their maximum value.
\end{teo}
The proof is simple, and can be reviewed at \cite{Gonzalez2018} it is based on the fact that the four triangles must be congruent to have the same energy.

Some measures $\mu(T)$ for triangles with those properties are

$$
\mu_1(T) = 4\sqrt{3} \frac{A}{l_1^2+l_2^2+l_3^2},  \quad \mu_2(T) = 2 \frac{r}{R},  \quad \mu_3(T)=\frac{4\sqrt{3}}{9} \frac{A}{R^2},  \quad\mu_4(T) = \frac{4}{\sqrt{3}} \frac{A}{l_{\max}^2}
$$
whree $\mu_1(T)$ is proposed by Joe, \cite{Joe2008},  $\mu_2(T)$ is the {\em radius ratio}  measure, $\mu_3(T)$ is  described by Shewchuk y $\mu_4(T)$ is the Cavendish's measure,  see \cite{Shewchuk2002}.

\begin{figure}[hbt]
\centerline{\begin{tabular}{cc}
\includegraphics[scale=.5]{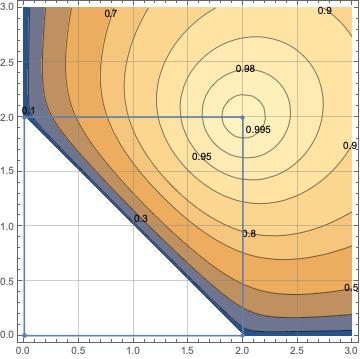} &
\includegraphics[scale=.5]{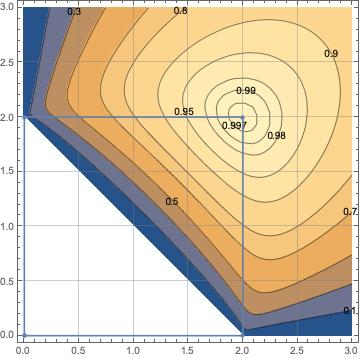} \\
(a) & (b)
\end{tabular}}
\caption{Level curves for (\ref{marmonica2v}) using (a) $\mu_1(T)$,  (b) {\em radius ratio} $\mu_2(T)$.}
\label{marmonica1}
\end{figure}

\section{Grid quality improvement}

\subsection{A distorsion measure}
In general, if $\mu(Q)$ is a good quality measure for quadrilaterals, a way of measuring the distortion of an $\epsilon$-convex quadrilateral $Q$ with respect to $\mu(Q)$ is using

$$
f(Q) = \frac{1}{\mu(Q)}
$$
because if $f(Q)$ is much greater than 1, the cell will be far from the value for which  $\mu(Q)$ characterizes the geometric shape of the cell $Q$ (square, rectangle, parallelogram, etc.) and we can say that $Q$ is a distorsioned quadrilateral with respect to that measure.

Under this idea, the distortion of the mesh $G$ can be measured as the average of the distortions of all the cells

\begin{equation}
F(G) = \frac{1}{Ne}\sum_{k=1}^{Ne} \frac{1}{\mu(Q_k)}
\label{e1c5}
\end{equation}
where $Ne$ its the number of the cells. We have in consequence the follow definition 

\begin{defi}
A grid $\hat{G}$ has better quality than the mesh $\bar{G}$ if

$$
F(\hat{G}) < F(\bar{G})
$$
where $F(G)$ is a distorsion measure.
\end{defi}

As an optimization problem,  improving the quality of a G-mesh can be considered as the problem

$$
G^* = \arg \min_G F(G) =  \frac{1}{Ne}\sum_{k=1}^{Ne} \frac{1}{\mu(Q_k)}.
$$
where the inner node of $G$ are the unknowns. The optimization problem is a large scale one when mesh dimension $m\times n$ is very large. It is important to notice that the initial mesh $G_0$ must be convex and remain so in each step of the optimization process.

Usually the quality measures for quadrilaterals are non differentiable functions, in an optimization process is better to build convex function with similar characteristics as the quality measure.

\section{New quality discrete functionals}
From the proof of theorem \ref{teoa} it is easy to see that
$$
2 a_c \leq \sqrt{(a^2+c^2)(b^2+d^2)} 
$$
for any convex quadrilateral; then

$$
f(Q)  = \frac{(a^2+c^2)(b^2+d^2)}{4 a_c^2}
$$
its a positive convex function whose critical points are rectangles. With this function we can define a discrete functional $F_R(G)$ over all  the grid cells

$$
F_R(G) = \frac{1}{Ne} \sum_{k=1}^{Ne}  f(Q_k).
$$

\begin{figure}[hbt]
\centerline{
\begin{tabular}{cc}
\scalebox{.5}{\includegraphics{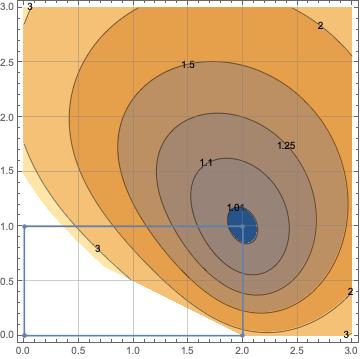}}
& 
\scalebox{.5}{\includegraphics{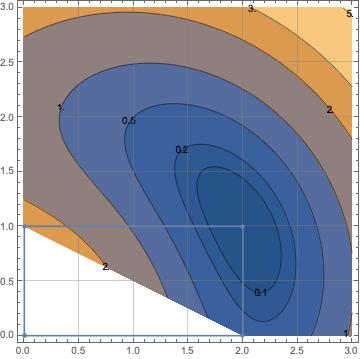}} \\
(a) & (b)
\end{tabular}}
\caption{Level for (a) $f(Q)$ y (b) $f_r(Q)$ where $Q$ has 3 fixed vertices $(0,0), (2,0)$ and $(0,1)$.} \label{nao}
\end{figure}

In Figure~\ref{nao} the shape of the surface of $F_R$ is sketched.

We propose to combine this functional with a convex area functional $S_w(G)$ (see \cite{Barrera2010}) to guarantee the convexity of the cells and and that their shape is close to a rectangle:

$$
F(G) = (1-\sigma) S_w(G) + \sigma F_R(G)
$$ 
where $\sigma>0$. In addition, the function 

$$
f(Q)  = \frac{(a^2+c^2)(b^2+d^2)}{4 a_c^2}
$$
can be interpreted (by cells) as a normalization (with respect to the jacobian) of Knupp's area-orthogonility functional \cite{Knupp1992}
$$
f_{ao}(Q)  = (a^2+c^2)(b^2+d^2).
$$

As we have discussed, because quality measures are usually non-differentiable functions it is difficult to use them as objective function; it is advisable to design convex and differentiable functions $f(Q)$ whose optimal values also satisfy $\mu(Q)\approx 1$ for a specific quality measure $\mu(Q)$.

\begin{figure}[hbt]
\centerline{\scalebox{.85}{\includegraphics{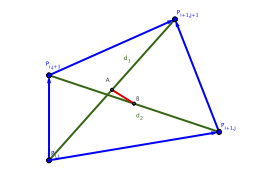}}}
\caption{Diagonals of a quadrilateral and the segment joining the midpoints between them.}
\label{diagonales5}
\end{figure}

As it is known, a rectangle is a parallelogram where its diagonals are equal and a quadrilateral is a parallelogram if its diagonals intersect at the midpoints. In Figure~\ref{diagonales5}, you can see those elements.

With this ideas, we propose a convex functional of the form

$$
F_p(G) = \sum_{i,j} \| \frac{{\bf p}_{i,j} + {\bf p}_{i+1,j+1}}{2}-\frac{{\bf p}_{i+1,j} + {\bf p}_{i,j+1}}{2} \|^2 
$$
over all elements of the grid $G$. Locally, this functional has as a critical point in the cells formed by parallelograms (see Khattri \cite{Khattri2006}). In that paper Kattri, uses a discrete scheme that coincides with the Laplacian-Isoparametric scheme of Herrmann \cite{Herrmann1976}.

Optimizing $F_p(G)$ tries to produce parallelograms, now we define a discrete functional to obtain rectangles. For each cell of $G$ let us measure the square of the difference of the square of diagonals

$$
F_d(G) = \sum_{i,j} (\| {\bf p}_{i,j} - {\bf p}_{i+1,j+1}\|^2- \| {\bf p}_{i+1,j} - {\bf p}_{i,j+1}\|^2 )^2;
$$
combining both functionals we get

$$
F_r(G) = \alpha F_p(G) + \beta F_d(G).
$$
If $\alpha\geq0, \beta \geq 0$ are appropiately chosen, we obtain a positive and convex functional  which has as a critical point in
a mesh formed by rectangles (including squares).  We can always do this if we guarantee that in each optimization step the mesh is convex.

Therefore, we use $S_w(G)$ to guarantee and preserve the convexity of the mesh, and combine with the latter fuctional to obtain a convex mesh and over irregular regions we will obtain meshes whose cells are close to being rectangles.

$$
F(G) = (1-\sigma)S_\omega (G)  + \sigma F_r(G).
$$

\newpage
\section{Examples}
\begin{figure}[hbt]
\centerline{\includegraphics[scale=.425]{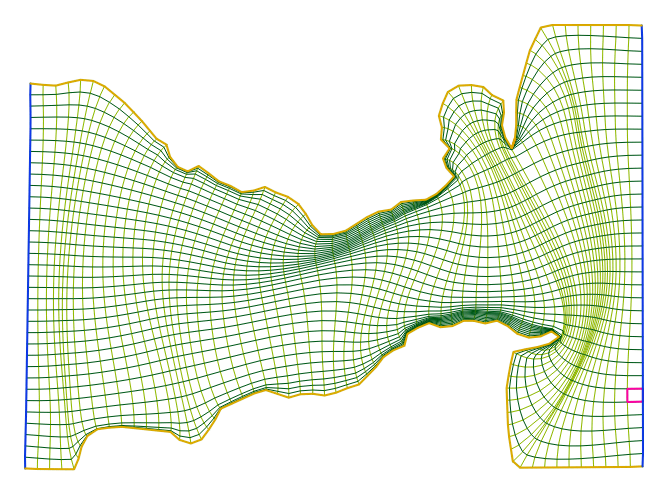}}
\caption{A mesh over the Strait of Gibraltar.}
\label{arlarge}
\end{figure}

For both functionals we can obtain optimal grids whose cells have a very large aspect ratio, see the Figure~\ref{arlarge}

To control the aspect ratio we propose use an area (volume) distortion measures functional
$$
F_A(G)=\sum_{q=1}^N  \frac{1}{\alpha(\triangle_q)^2} + \delta \alpha(\triangle_q)^2
$$
over all the $N$  signed areas of all the grid cells triangles.  Here $\delta>0$ is an adecuate value.  This distortion measures have a barrier on the boundary of the set of grids consisting of convex quadrilateral cells and its very similar to one proposed by Garanzha \cite{Garanzha2000},  but now we have a better control of the global distribution of the area. 

\begin{figure}[hbt]
\begin{center}
\begin{tabular}{c}
\includegraphics[scale=.225]{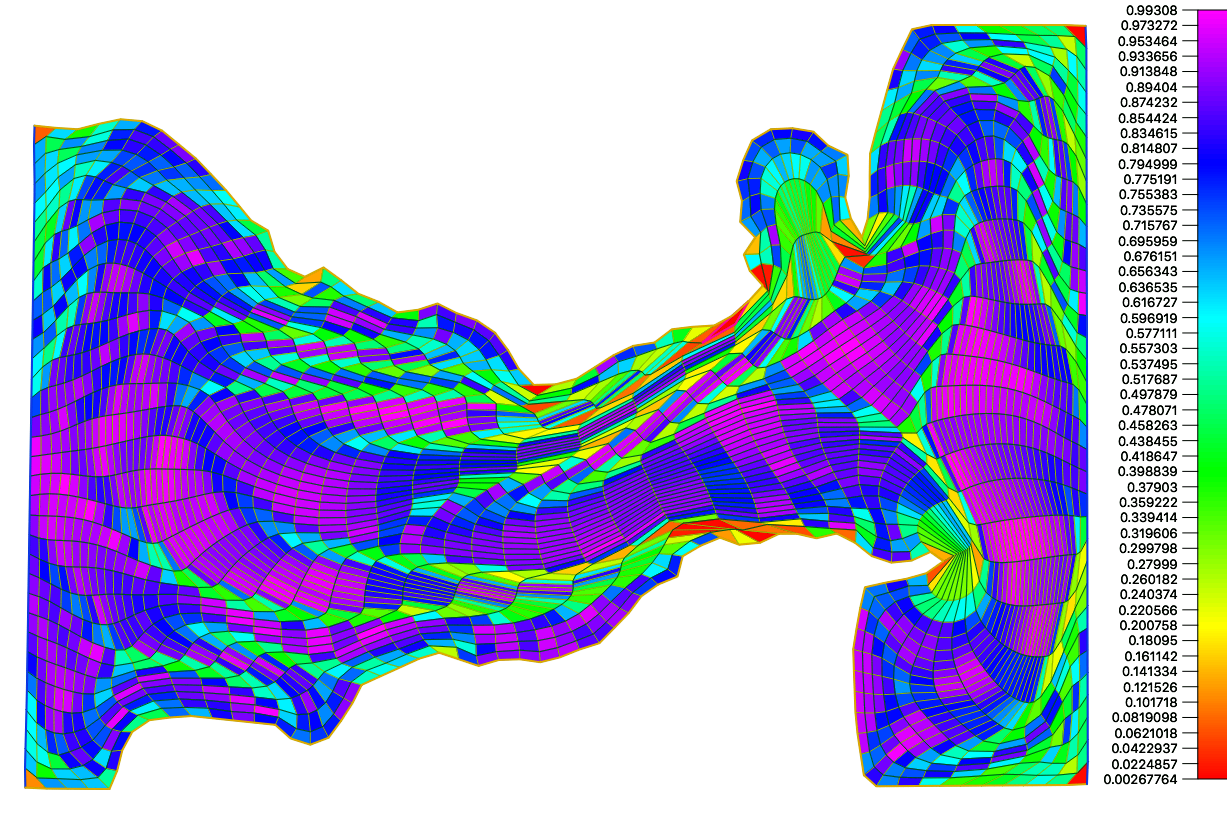} \\ (a) \\
\includegraphics[scale=.225]{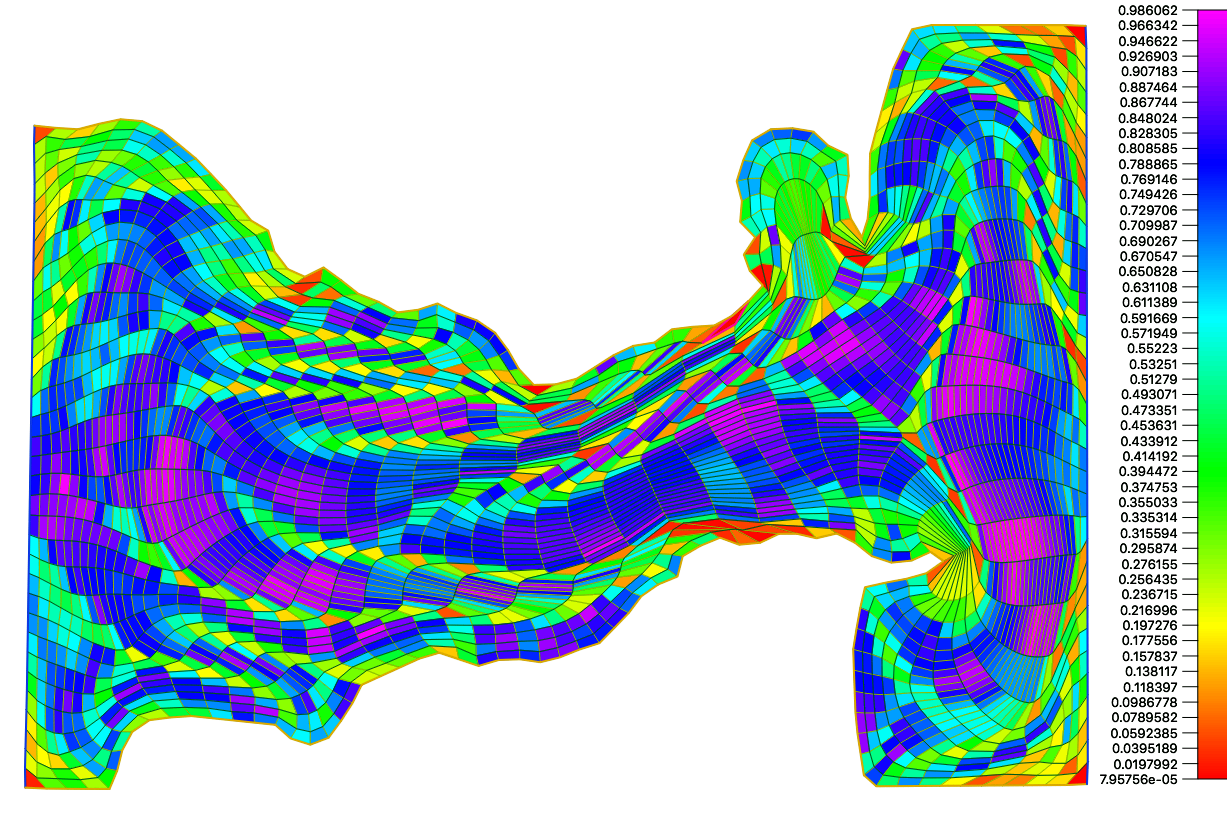} \\
 (b)
\end{tabular}
\end{center}
\caption{Color map of: (a) Rectangles 2015 quality measure (b) Minrect 2015 quality measure.}
\end{figure}

\begin{figure}[hbt]
\begin{center}
\begin{tabular}{c}
\includegraphics[scale=.3325]{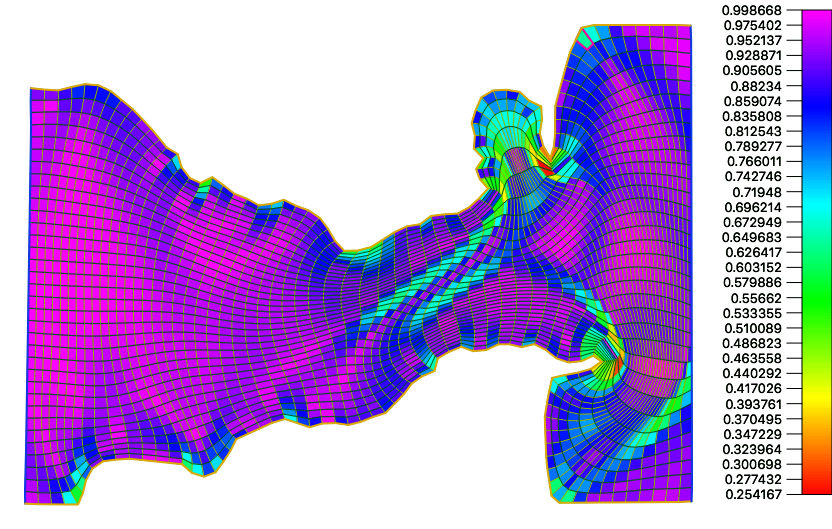} \\ (a) \\
\includegraphics[scale=.3325]{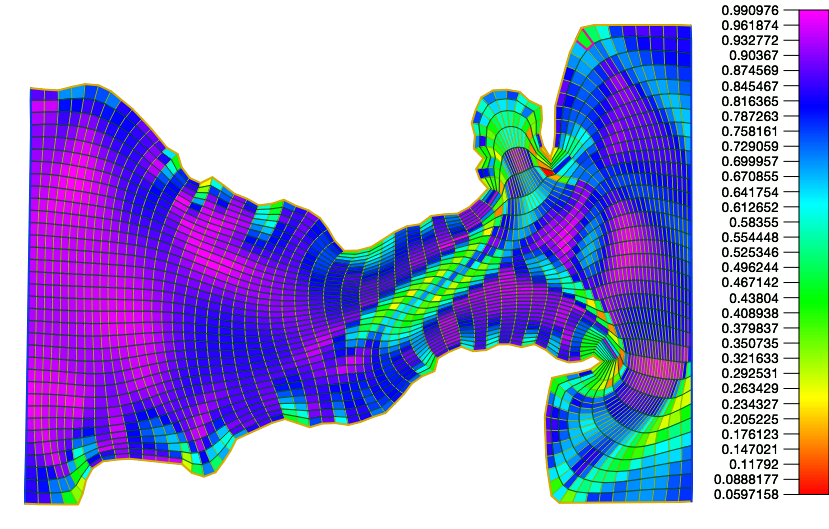} \\
 (b)
\end{tabular}
\end{center}
\caption{Color map of: (a) Rectangles 2015 quality measure (b) Minrect 2015 quality measure.}
\end{figure}

\section{Conclusions}
In this work we have described new quality measures for quadrilaterals, that help us to improve mesh and aspect ratio. Also we have proposed new functionals for grid generation that our computational tests, give us bases for using as alternative for area-orthogonal grid generation.

\section*{Acknowledgements} 
We want to thank Posgrado en Ciencias Matem\'aticas of UNAM for the financial support for this work. %We really appreciated the always valuable suggestions and comments made to this work by Dr. Fco. Dom\'\i nguez-Mota.

%\subsection*{Author contributions}
%
%This is an author contribution text. This is an author contribution text. This is an author contribution text. This is an author contribution text. This is an author contribution text. 
%
%\subsection*{Financial disclosure}
%
%None reported.

\subsection*{ORCID}
Guilmer Gonz\'alez Flores \url{https://orcid.org/0000-0001-7390-9088}
%
%\aiOrcid
%{\textcolor{orcidlogocol}{\aiOrcid} https://orcid.org/0000-0000-0000-0000}

\nocite{*}% Show all bib entries - both cited and uncited; comment this line to view only cited bib entries;
\bibliography{NQM-AMA}%

\end{document}